\documentclass{amsart}

\usepackage{amssymb}

\title{Options for a finite group model of quantum mechanics}
\author{Robert Arnott Wilson}
\address{Queen Mary University of London}
\email{r.a.wilson@qmul.ac.uk}

\date{First draft, 19th April 2021. This version, 7th July 2021.}

\begin{document}
\begin{abstract}
There are four finite groups that could plausibly play the role of the spin group in a finite or discrete model of quantum
mechanics, namely the four double covers of the three rotation groups of the Platonic solids. In an earlier paper I have
considered in detail how the smallest of these groups, namely the binary tetrahedral group,
of order $24$, could give rise to a non-relativistic theory that contains much of
the structure of the standard model of particle physics. In this paper I consider how one of the two double covers of the
rotation group of the cube might extend this to a relativistic theory.
\end{abstract}
\setcounter{tocdepth}{1}

\maketitle
\tableofcontents
\section{Introduction}
Most of theoretical physics, including quantum mechanics, quantum field theory \cite{Woit} and the
standard model of particle physics \cite{Griffiths}, is based on differential equations in continuous variables.
These equations naturally give rise to real and complex Lie algebras and Lie groups, and a great deal of geometry that is
built on top of these. But, as Einstein recognised \cite{Einstein1935} as early as 1935, if not earlier,
this creates a conflict with the fundamentally discrete nature of quantum mechanics. 
Nevertheless, no satisfactory discrete theory of quantum mechanics has been proposed, or at least accepted,
and the continuous theories continue to be used extremely successfully.

At its most basic group-theoretical heart, the theory is based on the rotation symmetry group $SO(3)$ of $3$-dimensional
Euclidean space, and its double cover $SU(2)$ that is required for modelling the spin $1/2$ properties of the
electron and other elementary particles. Any potential discrete version of this theory must therefore be based on a finite group
in place of $SO(3)$. 

In principle, any finite group could be considered, but there are really only three plausible candidates,
namely the three rotation symmetry groups of the five Platonic solids. The geometry of these solids is expounded in
detail in the 13th book of Euclid \cite{Euclid}, which dates back to around 2300 years ago.
The idea of abstracting these symmetries into the concept of a group did not take hold until the 19th century \cite{Galois},
after which there were serious studies of the tetrahedral group, of order $12$, the octahedral group, of order $24$,
and the icosahedral group, of order $60$, the last occupying the entirety of a famous book by Klein \cite{Klein}.

The octahedral group turns out to have two double covers, while the other two groups have one each, so that there are 
exactly four plausible options for a finite analogue of $SU(2)$ on which to try to build a finite model of quantum mechanics.
These groups have orders $24$, $48$, $48$ and $120$. In \cite{finite} I examined the smallest case in detail,
and found that much of the structure of the standard model of particle physics could be found in the representation
theory of this group. 
However, there were a number of issues to do with real versus complex representations,
and the fact that the Lorentz group did not seem to be in the right place, if indeed it was really there at all.

These problems suggest that an extension to
one of the groups of order $48$ might be required for a more satisfactory theory. In this paper I consider both
of these groups, and suggest that one of them is much more likely to work than the other. From this choice I then
show how it is possible to reconstruct most of the continuous structures that are used in the
standard theories. On the particle physics side this includes the Dirac spinors, the relativistic spin group, much of the
Dirac algebra, some version(s) of the Dirac equation and all of the gauge groups. On the relativity side, it includes the Lorentz group, 
which was problematic in the case of the smaller group, and possibly also the group of general covariance.

\section{Mathematical requisites}
\subsection{The groups}
All four of the suggested finite alternatives to $SU(2)$ contain the quaternion group, that consists of the $8$ unitary matrices of the form
\begin{eqnarray}
M(a,b)&:=& \begin{pmatrix}a & b\cr -\bar b& \bar a\end{pmatrix}
\end{eqnarray}
in which one of $a$, $b$ is $0$ and the other is $\pm 1$ or $\pm \mathrm i$.
For definiteness, let us define
\begin{eqnarray}
i:=M(0,1), & j:=M(0,\mathrm i), & k:=M(\mathrm i,0),
\end{eqnarray}
so that these matrices satisfy the familiar relations
\begin{eqnarray}
i^2=j^2=k^2=-1, && ij=-ji=k, jk=-kj=i, ki=-ik=j.
\end{eqnarray}

To obtain the double cover of the tetrahedral group, we may adjoin the matrix
\begin{eqnarray}
w&:=& (-1+i+j+k)/2\cr
&=& M(-1+\mathrm i, 1+\mathrm i)/2
\end{eqnarray}
so that we obtain in total $16$ new elements $M(\pm 1\pm \mathrm i, \pm 1\pm \mathrm i)/2$.
To obtain the two double covers of the octahedral group, we may adjoin one of the two matrices
\begin{eqnarray}
c&:=&(j-k)/\sqrt 2\cr
&=&M(-\mathrm i,\mathrm i)/\sqrt 2\cr
d&:=&c/\mathrm i.
\end{eqnarray}
For completeness, I give generators also for the double cover of the icosahedral group, although I shall
not use this group any further in this paper. Following \cite{icosians}, we may adjoin the matrix
\begin{eqnarray}
(i+\sigma j+\tau k)/2 &=& M(\tau \mathrm i, 1+\sigma\mathrm i)/2,
\end{eqnarray}
where 
\begin{eqnarray}
\sigma,\tau &:=& (1\pm \sqrt 5)/2.
\end{eqnarray}

The group that I shall be most concerned with in this paper is the only one of these four that does not lie inside $SU(2)$,
that is the group generated by abstract elements $i,w,d$ that satisfy the same
relations as the given matrices. The most important of these relations are as follows:
\begin{eqnarray}
iw=wj, & jw=wk, & kw=wi,\cr
id=-di, & jd=-dk, & kd=-dj,\cr
wdw=d, &w^3=1, & d^2=1.
\end{eqnarray}

\subsection{Conjugacy classes and representations}
The $48$ elements of the group 
\begin{eqnarray}
G&:=&\langle i,w,d\rangle
\end{eqnarray} 
fall into $8$ conjugacy classes, containing the following elements of the given orders:
\begin{eqnarray}
\begin{array}{ccc}
\mbox{Order} & \mbox{Elements} & \mbox{Number}\cr\hline
1 & 1 & 1\cr
2 & -1 & 1\cr
4 & \pm i,\pm j,\pm k & 6\cr
3 & w,iw,jw,kw,w^2,-iw^2,-jw^2,-kw^2 & 8\cr
6 & -w,-iw,-jw,-kw,-w^2,iw^2,jw^2,kw^2 & 8\cr
2 & \pm d, \pm id, \pm wd, \pm kwd, \pm w^2d, \pm jw^2d & 12\cr
8 & jd, -kd, iwd, -jwd, kw^2d, -iw^2d & 6\cr
8 & -jd, kd, -iwd,jwd,-kw^2d,iw^2d & 6\cr\hline
\end{array}
\end{eqnarray}

The representation theory \cite{JamesLiebeck}
of the group \cite{Zee} is summarised in the character table, in which I follow the usual mathematical convention
of indexing the columns by the conjugacy classes and the rows by the irreducible (complex) representations. The entries in the table
are the traces of the representing matrices.
\begin{eqnarray}
\begin{array}{c|cccccccc}
&1&-1&i&w&-w &d & jd & -jd\cr\hline
1^+&1&1&1&1&1&1&1&1\cr
1^-&1&1&1&1&1&-1&-1&-1\cr
2^0&2&2&2&-1&-1&0&0&0\cr
3^+&3&3&-1&0&0&1&-1&-1\cr
3^-&3&3&-1&0&0&-1&1&1\cr
2^+&2&-2&0&-1&1&0&\mathrm i\sqrt{2}&-\mathrm i\sqrt{2}\cr
2^-&2&-2&0&-1&1 &0&-\mathrm i\sqrt{2}&\mathrm i\sqrt{2}\cr
4^0&4&-4&0&1&-1&0&0&0\cr\hline
\end{array}
\end{eqnarray}
Only the last three representations are faithful, while the others represent various quotient groups. The representations $3^\pm$ represent the
octahedral group, isomorphic to the symmetric group 
$Sym(4)$; the representation $2^0$ represents $Sym(3)$, while $1^-$ represents $Sym(2)$ and $1^+$
represents the trivial group $Sym(1)$.

Since the group $G$ is generated by $w$ and $jd$, it is sufficient to specify the representing matrices for these two elements in each 
irreducible representation. The representations of dimensions $1$, $2$ and $3$ can be written as follows. 
\begin{eqnarray}
\begin{array}{c|cc}
&w & jd\cr\hline
1^\pm & (1) & \pm(1)\cr
2^\pm & M(-1+\mathrm i,1+\mathrm i)/2 & \pm\mathrm i M(\mathrm 1,\mathrm 1)/\sqrt 2\cr
2^0 & \frac12 \begin{pmatrix}-1&\sqrt 3\cr -\sqrt 3 & -1\end{pmatrix} & \begin{pmatrix}1&0\cr 0&-1\end{pmatrix}\cr
3^\pm &\begin{pmatrix}0&1&0\cr 0&0&1\cr 1&0&0\end{pmatrix} & \pm\begin{pmatrix}-1&0&0\cr0&0&-1\cr 0&1&0\end{pmatrix}\cr
\end{array}
\end{eqnarray}
Notice in particular that the representations $2^0$ and $2^\pm$ are unitary in the sense that in each case
the matrix group lies inside $U(2)$, but it does not
lie inside $SU(2)$. 
This 
seems to be an important property that is central to the construction of the electroweak gauge group
and the Dirac spinors. 

The representation $4^0$ can be constructed from $2^\pm$ by multiplying $w$ by a scalar $(-1\pm\mathrm i\sqrt3)/2$,
and multiplying $d$ by the complex conjugation map 
\begin{eqnarray}*&:&\mathrm i \mapsto -\mathrm i
\end{eqnarray}
and then expanding $\mathrm i$ to the $2\times 2$ matrix $i$, and expanding $*$ to the
obvious diagonal matrix. 
For a certain choice of signs, the matrices 
representing $w$ and $jd$ are then
\begin{eqnarray}
w&\mapsto&\frac14\begin{pmatrix}1-\sqrt3 & -1-\sqrt3 & -1-\sqrt3 & -1+\sqrt3\cr
1+\sqrt3 & 1-\sqrt3 & 1-\sqrt3 &-1-\sqrt3\cr
1-\sqrt3 &-1-\sqrt3&1+\sqrt3 & 1-\sqrt3\cr
1+\sqrt3 & 1-\sqrt3 & -1+\sqrt3& 1+\sqrt3\end{pmatrix} \cr
jd&\mapsto&
\frac1{\sqrt2}\begin{pmatrix}0&1&0&1\cr 1&0&1&0\cr 0&1&0&-1\cr 1&0&-1&0\end{pmatrix}
\end{eqnarray}

\subsection{A quaternionic notation}
One can see from the character table that there are four different maps from $G$ to $SO(4)$ that are of particular
significance, namely two faithful (i.e. fermionic) embeddings of $G$ with characters $2^++2^-$ and $4^0$, and two bosonic
representations (of $Sym(4)$) with characters $1^++3^-$ and $1^-+3^+$. One can also include the representation
of $Sym(3)$ with character $1^++1^-+2^0$. Since $SO(4)$ is generated by left and right
multiplications by quaternions, all five of these representations can be written in quaternionic notation. They are not quaternionic
representations, however, so that both left and right multiplications are required in all cases. 
The following table lists the images of a quaternion $q$ under the group generators $i,j,w,d$ in terms of quaternions
$i,j,k,w,c$ as defined above.
\begin{eqnarray}
\begin{array}{c|cccc}
&i&j&w&d\cr\hline
1^++1^-+2^0 & q & q & w^{-1}qw & c^{-1}qc\cr
1^++3^-& i^{-1}qi & j^{-1}qj & w^{-1}qw & c^{-1}qc\cr
1^-+3^+&i^{-1}qi & j^{-1}qj & w^{-1}qw & cqc\cr
4^0& qi&qj & w^{-1}qw & c^{-1}qc\cr
2^++2^-& qi & qj & qw & kqc \cr\hline
\end{array}
\end{eqnarray}

The basic boson/fermion distinction is the distinction between $-1$ acting as $+1$ or $-1$. In the context of the continuous
spin group $SU(2)$, this is the same as the distinction between elements acting by conjugation $q\mapsto x^{-1}qx$ or by
multiplication $q\mapsto qx$. But in the context of the finite group, the distinction is much more subtle. Of particular interest
is the fact that $4^0$ is `fermionic' for $i,j,k$, but `bosonic' for $w,d$: this fact will turn out to be very important for the physical
interpretation of this representation. 
If, for example, we want to interpret $4^0$ as a real version of the Dirac spinor, then $w$ and $d$ can be interpreted as acting on
`half' the spinor, since $w$ fixes the quaternions $1$ and $i+j+k$, while $d$ fixes the quaternions $1$ and $j-k$. But since $w$
and $d$ do not commute, the two `halves' are different. 
In effect, the standard model chooses the halves defined by $d$, and therefore has to abandon the symmetry $w$, which in
the group algebra model must represent generation symmetry.

Notice also the extra left-multiplication by $k$ in $2^++2^-$: this 
is another representation that we might want to interpret as a spinor, either a real Majorana spinor as here, or a complex Weyl spinor
$2^+$ or $2^-$, or even two Weyl spinors to make a Dirac spinor. 
But here $w$ and $d$ act without fixed points, so the only possibility for a consistent interpretation 
would seem to be as the left-handed
Weyl spinor.

The subtle differences between $2^\pm$ and $4^0$ as different types of spinors will play an important role in this paper.
Both of them contain important parts of the Dirac spinor, but neither of them individually can implement all aspects of the
Dirac spinor in the standard model. Nor is the distinction between them 
as simple as the distinction between left-handed and right-handed spinors.
Conversely, the Dirac spinor cannot implement all aspects of the fermionic representations of $G$.

The representation $1^++3^-$ is the finite analogue of the 
classical $4$-vector representation of $SU(2)$, distinguished from $1^-+3^+$ which plays the
role of the axial vectors.  The distinction between $1^++3^-$ and $1^++1^-+2^0$ is even more interesting, where replacing the
vector representation of $i,j,k$ by the scalar representation causes a symmetry-breaking of $3^-$ into $1^-+2^0$ which will play a vital
role in the implementation of electroweak symmetry-breaking and `mixing' of the 
forces.

\subsection{The structure of the group algebra}
From the character table we can read off the structure of the group algebra. The complex group algebra is simply a direct sum of complex
matrix algebras of the same sizes as the irreducible representations, that is
\begin{eqnarray}
&2\mathbb C + 3M_2(\mathbb C) + 2M_3(\mathbb C) + M_4(\mathbb C).
\end{eqnarray} 
But for physical applications it is likely that the real group algebra will be more useful. This has the structure
\begin{eqnarray}
&2\mathbb R + M_2(\mathbb R) + M_2(\mathbb C) + 2M_3(\mathbb R) + M_4(\mathbb R).
\end{eqnarray}
For mapping into the standard model, the two pairs of isomorphic algebras, of $1\times 1$ and $3\times 3$ real matrices, will need
to be converted to $1\times 1$ and $3\times 3$ complex matrices, in order to produce the gauge groups $U(1)$ and $SU(3)$
of electromagnetism and the strong force. Then $M_2(\mathbb R)$ is available for the group $SL(2,\mathbb R)$ as a broken symmetry version
of the weak gauge group, while $M_2(\mathbb C)$ produces the relativistic spin group $SL(2,\mathbb C)$. 
This leaves $M_4(\mathbb R)$ for a real version of the Dirac algebra, if we can find a way of implementing the complex
structure somewhere else.
If so, then this provides all the
necessary ingredients for the standard model of particle physics. 

Now we must be careful to distinguish two completely different types of multiplication.
The first is the group multiplication, which extends linearly to the whole group algebra, and is always implemented as
matrix multiplication, in every representation. The second is the tensor product of representations. The Dirac algebra, being
a Clifford algebra, includes some aspects of this second multiplication as well, also implemented as matrix multiplication.
But in general, tensor products of representations cannot be implemented via matrix multiplication. This is the reason why, for example,
the strong force cannot be implemented inside the Dirac algebra.

If we want to implement the Dirac algebra in terms of the group multiplication, then it must essentially be the matrix algebra
$M_4(\mathbb R)$, with an extra complex structure arising from somewhere to be determined. If on the other hand it is the
structure of the algebra as the tensor product of two copies of the Dirac spinor representation that is important, then we must 
instead take the
tensor product of the finite group representations
\begin{eqnarray}
2^+\otimes 2^+\cong 2^-\otimes 2^- &=& 1 ^-+3^+\cr
2^+\otimes 2^-\cong 2^-\otimes 2^+ &=& 1^++3^-.
\end{eqnarray}
Or, if we actually need aspects of both versions, then we also need to consider the tensor square of $4^0$, that is
\begin{eqnarray}
4^0\otimes 4^0 &=& (1^++3^-)+(1^-+3^+) + (2^0 + 3^++3^-). 
\end{eqnarray}
Notice that $4^0\otimes 4^0$ differs from $(2^++2^-)\otimes(2^++2^-)$ only in the replacement of $1^++1^-$ by $2^0$.

\subsection{Quotient groups and subalgebras}
As already noted, the structure of the group algebra as a direct sum of subalgebras
is closely related to the structure of the group and its normal subgroups, and corresponding quotient groups.
The quotient groups form a single series
\begin{eqnarray}
&G \rightarrow Sym(4) \rightarrow Sym(3) \rightarrow Sym(2) \rightarrow Sym(1)
\end{eqnarray}
with corresponding normal subgroups
\begin{eqnarray}
&1 \subset Z_2 \subset Q_8 \subset 2.Alt(4) \subset G.
\end{eqnarray}

The group algebra of each {quotient} group is both a quotient algebra and a subalgebra of the group algebra of $G$. 
As we proceed down the chain of quotients, therefore, we lose more and
more of the algebra:
\begin{eqnarray}
&&2\mathbb R + M_2(\mathbb R) + 2M_3(\mathbb R) + M_2(\mathbb C) + M_4(\mathbb R)\cr
&\rightarrow&2\mathbb R + M_2(\mathbb R) + 2M_3(\mathbb R) \cr
&\rightarrow&2\mathbb R + M_2(\mathbb R)\cr
&\rightarrow&2\mathbb R \cr
&\rightarrow&\mathbb R 
\end{eqnarray}
The group algebras of the corresponding normal subgroups are much less closely related to the group
algebra of $G$, and play very little role in the theory that I aim to describe. 
For reference, we have
\begin{eqnarray}
\mathbb R 1 &=& \mathbb R\cr
\mathbb  R Z_2 &=& 2\mathbb R\cr
\mathbb R Q_8 &=& 4\mathbb R + \mathbb H\cr
\mathbb R 2.Alt(4) &=& \mathbb R +\mathbb C + M_3(\mathbb R) + \mathbb H + M_2(\mathbb C).
\end{eqnarray}

The three nontrivial proper quotient groups $Sym(2)$, $Sym(3)$ and $Sym(4)$ appear therefore to be closely related to the
 subalgebras $2\mathbb R$, $M_2(\mathbb R)$ and $2M_3(\mathbb R)$ respectively. In the standard model, these
 subalgebras appear in complex form as
 \begin{eqnarray}
 \mathbb C + \mathbb H + M_3(\mathbb C), 
 \end{eqnarray}
 with the same real dimensions, but with an imposed unitary structure. The original purpose of this complex structure seems to
 have been to enable the quotient group $Sym(2)$ to be implemented in electromagnetism as complex conjugation. This was then
 abstracted to a unitary group $U(1)$ in order to provide something for the finite group to act on. By analogy, the other
 algebras give rise to unitary subgroups $SU(2)$ and $U(3)$. Hence we obtain a `gauge group'
 \begin{eqnarray}
 U(1) \times SU(2) \times U(3)
 \end{eqnarray}
 which describes much of the structure of the bosonic part of the algebra, although it only contains $13$ of the total of $24$ dimensions.
 
 This part of the standard model works well. But difficulties arise when these Lie groups are called upon to play the role of the
 finite groups as well. Thus the Lie group $SU(2)$, corresponding in the group algebra to the finite group $Sym(3)$, is called
 upon to define `weak doublets' and play the role of the group $Sym(2)$. This is, I believe, the precise point at which the three-generation
 structure is lost, and the three copies of $Sym(2)$ inside $Sym(3)$ become identified with each other. More precisely, there is an
 ambiguity as to whether $Sym(2)$ should be considered as a subgroup  of $Sym(3)$ or as a 
 quotient group of $Sym(3)$, or perhaps both.
 
 A similar but more complicated situation arises when the Lie group $U(3)$, or $SU(3)$, is called upon to act as $Sym(3)$
 on the three generations, rather than as $Sym(4)$. The full implications of this cannot be considered until we have resolved the
 problem of the mixing of the quotients $Sym(3)$ and $Sym(2)$ satisfactorily.
 
 \section{The standard model}
 \subsection{Mixing}
 There is therefore a kind of `mixing' between the various gauge groups arising from the action of
the finite group. To describe this, it is not necessary to use the whole of the group algebra, but only the ambient orthogonal
groups that act on the irreducible representations. Omitting the trivial representation, where there is no mixing to consider, we have the groups
\begin{eqnarray}
O(1), O(2), O(3), SO(3)
\end{eqnarray} 
 acting on $1^-$, $2^0$, $3^+$ and $3^-$ respectively.
 The relationship between $O(1)$ and $O(2)$ is mathematically most naturally expressed as a quotient map from $O(2)$ to $O(1)$, defined by the
 determinant. 
 
  In the standard model, however, the corresponding relationship between $U(1)$ and $U(2)$ is expressed in terms of a subgroup.
 This involves a choice of a particular reflection in $O(2)$ to identify with $O(1)$, and this choice can be expressed in terms of
 an angle, in this case the Weinberg angle. 
 But this choice depends on a choice of $d$ as one of the three involutions in the finite group $Sym(3)$,
 and therefore on a choice of one of the three generations of electrons. Hence the choice depends in some way
 on mass, and therefore on energy. So it is not entirely surprising to see an experimental dependence of the
 Weinberg angle on the energy involved in any given experiment.
 
 Similarly, the relationship between $O(2)$ and $O(3)$ is expressed as a subgroup relationship. One needs two angles to specify the
 direction of the rotation subgroup $SO(2)$, and a third to specify the angle of rotation. 
 To include the reflections, one needs to specify one more angle, between the actions of $d$  
as reflections in $O(2)$ and  $O(3)$. At this point it is worth noting that $d$ acts in $O(1)$, and therefore in $O(2)$, as charge conjugation,
whereas in $O(3)$ it acts to reverse the parity. Therefore the angle between these two actions of $d$ is the angle between the
C and P symmetries in the standard model, known as the CP-violating phase. 

The other three angles are most likely best
identified as the remaining three angles in the 
Cabibbo--Kobayashi--Maskawa (CKM) matrix, of which the most prominent is the Cabibbo angle \cite{Cabibbo,KM},
that describes `mixing' between first and second generation quarks. Using orthogonal groups rather than unitary groups 
explains why this unitary matrix
contains only four independent real parameters, rather than the nine that might be expected.
There is a further mixing between $O(3)$ and $SO(3)$ that might perhaps correspond to the 
Pontecorvo--Maki--Nakagawa--Sakata (PMNS) matrix \cite{Pontecorvo,MNS}, that attempts to quantify the
phenomena of neutrino oscillation in a similar way that the CKM matrix describes quark mixing.
 
\subsection{A subgroup}
Much of the structure of the standard model, such as the Dirac algebra of complex dimension $16$, only applies to one generation of
fermions at a time. In order to reproduce such structures from the finite group, it seems likely that we will need to take out the finite
symmetries of order $3$, and restrict therefore to a subgroup of order $16$, of which there are three. These subgroups and their group
algebras may give rise to some structures like the Dirac algebra. However, there must be some differences, since the Dirac algebra
as it stands arises as a twisted group algebra of an abelian group of order $16$, whereas our group is non-abelian. We can only hope
that this non-abelian structure can be used to implement some of the symmetry-breaking that otherwise has to be implemented by
hand. 

Let us take the subgroup $H$ generated by $j$ and $d$.  This group has order $16$, and its elements fall into the following seven conjugacy classes:
\begin{eqnarray}
\begin{array}{ccc}
\mbox{Order} & \mbox{Elements} & \mbox{Number}\cr\hline
1 & 1 & 1\cr
2 & -1 & 1\cr
4 & \pm i & 2\cr
4 & \pm j, \pm k & 4\cr
2 & \pm d, \pm i d& 4\cr
8 & jd, kd & 2\cr
8 & -jd, -kd & 2\cr\hline
\end{array}
\end{eqnarray}
This group is known as the semi-dihedral group of order $16$, since it can be generated by the elements $jd$ and $d$ that
satisfy the relations
\begin{eqnarray}
(jd)^8=1, &d^2=1,& (jd)d = d(jd)^3.
\end{eqnarray}
The quotient of this group by the group of signs is isomorphic to the dihedral group $D_8$ of order $8$, which gives $5$
non-faithful irreducible representations. There are then two faithful complex representations of dimension $2$, which we might hope
to identify with the Dirac spinor representation in some way.

The character table is as follows:
\begin{eqnarray}
\begin{array}{c|ccccccc}
&1&-1&i&j&d&jd&-jd\cr\hline
1a&1&1&1&1&1&1&1\cr
1b&1&1&1&1&-1&-1&-1\cr
1c&1&1&1&-1&1&-1&-1\cr
1d&1&1&1&-1&-1&1&1\cr
2a&2&2&-2&0&0&0&0\cr
2b&2&-2&0&0&0&\mathrm i\sqrt 2&-\mathrm i\sqrt 2\cr
2c&2&-2&0&0&0&-\mathrm i\sqrt2&\mathrm i\sqrt 2\cr\hline
\end{array}
\end{eqnarray}
Before we go further in the representation theory, however, it is worth looking more closely at  
the abstract structure of $H$.

\subsection{Dirac matrices}
Since $H$ acts on some representations that we want to identify with Dirac spinors, we need to investigate how $H$
is related to the finite (abstract) group 
generated by the Dirac matrices. Taking the latter to be
\begin{eqnarray}
\gamma_1,\gamma_2,\gamma_3,\gamma_0,&\gamma_5=\mathrm i \gamma_1\gamma_2\gamma_0\gamma_3
\end{eqnarray}
the important relations between them are that the squares are
\begin{eqnarray}
\gamma_1^2=\gamma_2^2=\gamma_3^2=-1,&\gamma_0^2=\gamma_5^2=1,
\end{eqnarray}
and any two anti-commute.
These relations imply that the finite group they generate has order $64$.

Now we have a group of order $16$, so we cannot find elements satisfying exactly these relations.
But we can come close, if we take the elements
\begin{eqnarray}
i,j,k,d,&id=ijdk.
\end{eqnarray}
The squares are the same, and the anti-commutation rules are the same, apart from the fact that
\begin{eqnarray}
jd=-dk,\quad kd=-dj,&&j(id)=(id)k,\quad k(id)=(id)j,
\end{eqnarray}
compared to
\begin{eqnarray}
\gamma_2\gamma_0=-\gamma_0\gamma_2,\quad
\gamma_3\gamma_0=-\gamma_0\gamma_3,&&
\gamma_2\gamma_5=-\gamma_5\gamma_2,\quad
\gamma_3\gamma_5=-\gamma_5\gamma_3.
\end{eqnarray}
These differences are actually very subtle, having to do with properties of changes in the direction of spin or momentum in
elementary particle interactions, particularly electro-weak interactions. When we put the continuous variables back in,
therefore, we might expect these differences to appear in the modelling of photon polarisation, and in properties of
spin-flipping in weak interactions, for example. Of course, we must be alert to the possibility that these changes
are detrimental to the model, and break some essential part of the standard model.

But there are at least three potential advantages, that suggest it is not time to give up yet:
\begin{itemize}
\item
we have only $Q_8$ generated by $i,j,k$, rather than $Z_2\times Q_8$
generated by the corresponding Dirac matrices;
\item we do not have to complexify the algebra in order to enforce the condition $\gamma_5^2=1$;
\item we have three copies of the group $H$ available for three generations, by choosing which one of $i,j,k$
to invert, and which two to swap.
\end{itemize}
The first two are technical simplifications which are of mainly mathematical rather than physical interest.
Between them, they reduce the size of the group by a factor of $4$, which could significantly simplify the calculations,
provided that they have not at the same time thrown away some vital piece of information.
The third, on the other hand, 
is potentially of huge 
importance for physics, because it incorporates the three generations into the
structure of the algebra in a natural way, so that, if it all works out, then it will no longer be necessary to do this by hand afterwards.

It is to be hoped that the above suggestions will be sufficient to implement the vertex factors of the Feynman calculus, once
we have sorted out how to extend from the finite group to the Lie groups, in order to implement the Weinberg angle. The propagators,
however, cannot be implemented just by matrix multiplication, but require some tensor product structure  
to produce the correct representations. 
In the group algebra of $H$, the tensors we might require for bosons are
\begin{eqnarray}
2b\otimes 2b \cong 2c\otimes 2c &=& 1b+1c+2a,\cr
2b\otimes 2c\cong 2c\otimes 2b &=& 1a+1d+2a.
\end{eqnarray}
These correspond to $1^\mp+3^\pm$, respectively, in the group algebra of $G$. 

In order to try to match this up with the even part of the Dirac algebra, we need an 
interpretation of the
elements of $H$. The interpretation of $i,j,k$ must surely be 
\begin{eqnarray}
\gamma_2\gamma_3 &\leftrightarrow& jk=i,\cr
\gamma_3\gamma_1 &\leftrightarrow& ki=j,\cr
\gamma_1\gamma_2 &\leftrightarrow& ij=k.
\end{eqnarray} 
It would now seem to be reasonable to interpret $id$ as $\gamma_5$, just as before.
Of course, in the standard model, $\gamma_5$ commutes with the even part of the Clifford algebra, whereas $id$
has a non-trivial action on $i,j,k$. 
This means that the finite group implements some discrete actions on certain quantum
numbers of the particles. 

This is, of course, the real motivation for studying finite groups here,
so that these discrete actions can be incorporated into the basic algebra, rather than having to be
added in by hand later.
If we can match up the discrete structure with the continuous structure
correctly, then we 
might be able to throw some light on the mechanism of electro-weak mixing.

\subsection{Symmetry-breaking}
The information we need for understanding how the symmetry-breaking works is the way the irreducible representations of $G$
restrict to $H$: 
\begin{eqnarray}
\begin{array}{c|cc|ccc|cc}
&1a&1b&1c&1d&2a&2b&2c\cr\hline
1^+ & 1&&&&&&\cr
1^- & & 1&&&&&\cr
2^0& 1&1&&&&\cr\hline
3^+&&&1&&1\cr
3^-&&&&1&1\cr\hline
2^+ &&&&&&1\cr
2^-&&&&&&&1\cr
4^0&&&&&&1&1\cr\hline
\end{array}
\end{eqnarray}
There are three 
separate `blocks' of symmetry-breaking, 
that look rather similar to each other.
The first involves the representations $1^\pm,2^0$ of $G$ and $1a/b$ of $H$, and we might expect it describe symmetry-breaking
of the weak force, together with electro-weak unification. The second involves the representations $3^\pm$ of $G$ and $1c/d,2a$ of $H$,
and might be expected to relate to the mixing of the electroweak and strong forces. The third involves the representations $2^\pm, 4^0$
of $G$, and $2b/c$ of $H$, so might relate the Dirac spinors to properties of gravity and/or mass.

It is also striking that the three `blocks' have much the same structure. The first block has the structure $1+1+2\times 2$ for $G$,
compared to three copies of $1a+1b$ for $H$. The third block has the same shape, but has been tensored with a Weyl spinor
to give a structure of three Dirac spinors for $H$, or one Dirac spinor plus a $4\times 4$ matrix algebra for $G$. The second block
has been transposed, so that the $1+1+2\times 2$ structure appears three times for $H$, rather than once for $G$, while for
$G$ we have a pair of $3\times 3$ algebras. 

In other words, there is a curious duality between the restriction maps on
the first and third blocks, and the induction map (as it is known) on the second block.
It is therefore unlikely to be possible to separate these three blocks from each other. But the first block looks like three generations of
a weak doublet of leptons, and after tensoring with a Weyl spinor gives the three Dirac spinors that describe these particles in the
standard model. The second block similarly looks like three generations of a weak doublet of quarks, with six copies of a
2-dimensional representation to hold the `colours'. In this interpretation, colour confinement is enforced by the fact that these representations
are $2$-dimensional, not $3$-dimensional.

Tensoring with a Weyl spinor obliterates all differences on the classes represented by $i$, $j$ and $d$, so that the
representations $1a+1b$ and $1c+1d$ and $2a$ are all mapped to Dirac spinors in the same way. In other words,
the Dirac spinor cannot distinguish between the three generations, cannot distinguish between leptons and quarks,
and cannot distinguish between the three colours. This is all consistent with the properties of the Dirac spinor in the standard model.
For an initial investigation, therefore, it may make sense to divide the group algebra of $G$ into three subalgebras:
\begin{enumerate}
\item the lepton subalgebra, isomorphic to the group algebra of $Sym(3)$;
\item the 
quark subalgebra, that extends to the group algebra of $Sym(4)$;
\item the 
spinor subalgebra, that extends to the full group algebra of $G$.
\end{enumerate}

\section{A toy model}
\subsection{
The lepton algebra}
We must now begin to grapple with the serious problem of relating the discrete to the continuous.
The group algebra is the natural structure in which such a relationship takes place, since it is a real algebra
with a finite group as a preferred basis. In representation theory the usual convention is for the finite group
to act by right multiplication on the algebra, and for the algebra (or any associated subalgebras or Lie groups)
to act on itself by left multiplication. In that way we can separate the `gauge groups'---really gauge algebras---acting on the left, from
the discrete symmetries acting on the right.

To see how this process works in detail, it is best to start with a toy model. For this purpose, the lepton
subalgebra is just big enough to exhibit the salient features, while being small enough to write down
all the calculations in full. This algebra is isomorphic to the group algebra of the group $K:=Sym(3)$, which we can take to
be generated by $w$ and $d$, as either a subgroup or a quotient group of $G$. The group consists of the
identity element $e=1$, two rotations, $w$ and $w^2=w^{-1}$, and three reflections, $d$, $wd=dw^{-1}$ and
$w^{-1}d=dw$.
There are three irreducible representations, namely the trivial representation $1^+$, the alternating representation $1^-$
and the deleted permutation representation $2^0$.

 Each representation defines an idempotent in the
group algebra, obtained by multiplying each group element by its trace in that representation, averaging over the group
and multiplying by the dimension:
\begin{eqnarray}
p&:=& (1+w+w^2+d+wd+w^2d)/6\cr
q&:=& (1+w+w^2-d-wd-w^2d)/6\cr
r&:=&(2-w-w^2)/3
\end{eqnarray}
These idempotents satisfy the relations
\begin{eqnarray}
p^2=p, q^2=q, r^2=r,& pq=qp=pr=rp=qr=rq=0,& p+q+r=1.
\end{eqnarray}
The group algebra $\mathbb R K$ then decomposes as the direct sum of the three subalgebras
\begin{eqnarray}
\mathbb Rp \cong \mathbb R,& \mathbb Rq \cong \mathbb R, & \mathbb RrK \cong M_2(\mathbb R).
\end{eqnarray}
The last of these isomorphisms can be defined on two generators of the algebra via 
\begin{eqnarray}
rd&\mapsto & \begin{pmatrix}1&0\cr 0&-1\end{pmatrix},\cr
r(w-w^2)/\sqrt3&\mapsto &\begin{pmatrix}0&1\cr-1&0\end{pmatrix},
\end{eqnarray}
after which the relations are easily verified.

The relations
\begin{eqnarray}
pw=pd=p, & qw=q, & qd=-q
\end{eqnarray}
give the actions of the finite group on the $1$-dimensional representations.
There are two copies of the $2$-dimensional representation, 
on the two rows of the matrix algebra. There is then a gauge group $SL(2,\mathbb R)$ that
acts by left-multiplication to mix these two copies together in arbitrary ways.
This gauge group has no effect on the action of the finite group, but may incorporate some
continuous variables that arise from the observer's choice of coordinate system.

To get a complete description of the algebra, therefore, we note that $d$ acts on $\mathbb Rp+\mathbb Rq$ in exactly the
same way that it acts on each of the rows of $\mathbb RrK$, while $w$ acts as $1$, so that $(w-w^2)/\sqrt3$ acts as $0$.  
This contrasts with the action as $1$ on the representation $2^0$, and when transferred to the spinor representation 
may be capable of interpretation as an action on the `left-handed' spinor only.

\subsection{Finite versus infinite}
The job of the finite group is to describe the discrete symmetries, that is the division of the leptons into three generations
of electrons and neutrinos. This division is defined by electro-weak interactions. The job of the gauge groups is to define
the coordinates that are used by a particular observer of a particular experiment in particular circumstances. 
The standard model uses the same Lie group for both purposes.

Essentially, the standard model starts with a gauge group $U(1)\times SU(2)$ acting on the
representation $1^++1^-+2^0$, and then 
multiplies $\gamma_1\gamma_2\gamma_0\gamma_3$ by the
complex number $\mathrm i$ in order to create the element $\gamma_5$ which is here implemented as $id$. 
Multiplication by $\mathrm i$ has the
effect of rotating the $2$-space through $90^\circ$. Then
in order to ensure that $\gamma_5$ commutes with the generation symmetry $w$, the standard model requires a further $30^\circ$
rotation to move $w$ to $1$. There is then a further correction to this $30^\circ$ angle to take account of the fact that the
symmetry between proton and neutron is not exact. As a result, the standard model ends up with a Weinberg angle of
approximately $29^\circ$, that varies with the energy scale.

This process has 
two notable 
effects on the finite group. First of all, by enforcing the relation $w=1$, the standard model has factored out the generation
symmetry, so that it is not able to incorporate three generations in the algebra. Second, when the extension to quarks requires $w$
to be reinstated, it becomes apparent that $id$ (or $\gamma_5$) does not in fact commute with the generation symmetry $w$,
and a Cabibbo angle has to be introduced to deal with this fact. At the lepton level, however, the only
effect is to swap the names of the $\mu$ and $\tau$ neutrinos. This difference may be important when considering 
neutrino oscillations \cite{oscillation,neutrinos,SNO}.

Returning now to the gauge group proper, or more precisely the gauge algebra $2\mathbb R + M_2(\mathbb R)$,  there are
three real scale factors, 
two of which define the units of mass and charge, so are of no physical interest. 
Then there is an electromagnetic dimensionless parameter that defines the unification of $\mathbb Rp+\mathbb Rq$
into a copy of the complex numbers, $\mathbb C$. This is presumably related in some way to the fine-structure constant $\alpha$, 
and the fact that $\alpha\ne0$ gives the electromagnetic $U(1)$ a slightly elliptical shape.
This leaves three more electro-weak parameters, which we might take as the two mass ratios between the three charged leptons,
together with the Weinberg angle.
Certainly there are other ways of choosing the $6$ fundamental parameters, but this seems a reasonable choice for this toy
model. In particular, the neutrinos in this model are massless, though 
they may acquire mass from somewhere else in the
full model. 

In this toy model it seems that at least some of these parameters can be taken as arbitrary
and that the theory is exactly the same whatever values are taken. 
This should not be taken to mean that the parameters necessarily \emph{do} vary between observers,
in any real sense, only that the theory may have more symmetry than is immediately apparent. 
If this symmetry survives into the full model, then it is a kind of extension to the general principle of relativity, of the same general form, that
although we measure particular space/time/mass/charge coordinates, the theory does not depend on our
choice of coordinate system. 

\subsection{Fundamental parameters}
If this is the case, then it begs the question, how are the actual measured values determined, and can we calculate them
from the model? In principle, the group algebra is simply a weighted average over the whole universe
of copies of the group elements, each copy
representing a single interaction between two leptons. 
Since the model includes massless neutrinos, which therefore travel at the speed of light, we must include quite a large part of
the universe in this averaging process. 

Certainly we must go at least to the Solar System scale, since most of the neutrinos
we can detect come from the Sun. Possibly we need to go to the scale of the galaxy, and maybe larger, perhaps
the entire universe, to take account of neutrinos emanating from the Big Bang.
Of course, it is impossible to isolate the leptons from all the other particles in the universe, so we cannot contemplate any such
calculations until we have a complete model, including a quantum theory of gravity. Even then, the calculations are likely to be
wildly impractical.

In any case, the group $SL(2,\mathbb R)$ that I am using here to gauge the leptons plays quite a different role in the
group algebra from the role it plays in the standard model. It may not therefore be possible to identify the two in any
meaningful way. The standard model $SL(2,\mathbb R)$ is used to define the three weak doublets, of electron and associated neutrino.
It is therefore playing the role of the finite group $Sym(2)$ acting on the right, rather than the Lie group acting on the left.

One can, of course, make a continuous group act on the right, by the process of quantum superposition. Or one can simply
interpret such a group as describing the average behaviour of a collection of similar particles in similar circumstances.
Using the group in this way allows one to calculate probabilities of certain outcomes, given a suitable probability distribution of
inputs. It also allows us to scale down the gauge group acting on the left, so that we no longer have to average over the whole
universe, but only over the scale of a particular experiment. 

But by doing so, we ignore the influence of the outside world,
which not only introduces some uncertainty into the calculations, but also introduces a systematic bias, which is corrected for
by the standard model parameters. Conversely, the greater the precision that is required for standard model parameters,
the larger the region of the universe that is required for averaging in order to obtain that precision. The standard model assumes that
this averaging can always be done locally, but unfortunately experiment long ago reached the point where that is no longer true.

\subsection{A possible application}
Let us take an example that has been in the news recently, and can be discussed purely in terms of the lepton algebra, that is the
muon $g-2$ experiment \cite{muong-2}. 
All the local influences have been meticulously accounted for \cite{muontheory,muonHVP}, and yet the experimental value of $g-2$ for the muon
is measured at around $(2.15\pm.35)$ppm higher than the calculated value. 
At least, this is the current consensus, although not everyone agrees.  There are competing calculations using lattice QCD
\cite{Fodor} which find no tension
between theory and experiment. 
Indeed, it may be that the discrete nature of lattice QCD gives a better approximation to the underlying reality.
For present purposes, I assume that the consensus is correct, though clearly it may not be.

In other words, all the known fundamental forces except gravity have been
fully accounted for. Of course, we haven't got a quantum theory of gravity yet, so we can't account for gravity in the same way.
But that does not mean that we have to ignore gravity, since we can 
use Newtonian gravity instead.
The muons 
in these experiments were in a horizontal orbit with a radius of approximately 15 meters.
The radius of the Earth varies between 6357 and 6378 km. The angle by which the direction of the gravitational field
changes over a distance of 15 meters is therefore around $2.35\times 10^{-6}$ radians. 

This numerical coincidence shows
that
 the experimentally measured effect 
is consistent with a gravitational origin of some kind. 
Of course, there is no theory of quantum gravity that can predict this coincidence, but 
before claiming any evidence for a fifth force, it is first necessary to rule out the fourth force. 
At the present time, such ruling out is 
both theoretically and practically impossible.

In this context it may be worth noting that a similar numerical coincidence occurs in the experiment that detected CP-violation
of neutral kaon decays \cite{CPexp}. In this case, the length of the experiment was reported as $57$ feet, so that the 
corresponding angle is around $2.73\times 10^{-6}$ radians. This is consistent with the reported proportion of the kaons
that appeared to have changed from the long eigenstate to the short eigenstate over that distance. Of course, there is no good
reason to think this is anything other than a meaningless coincidence, but perhaps there is some connection to other anomalies
for kaons \cite{kaonanomaly,Kaon2} and B-mesons.

\section{A bigger model}
\subsection{The quark algebra}
The next stage, then, is to examine the group algebra of $Sym(4)$. This adds $18$ dimensions, consisting of two copies of
the $3\times 3$ real matrix algebra, to the lepton algebra already described. In the standard model, this is converted into the
$3\times 3$ complex matrix algebra in the same way that $2\mathbb R$ was converted into $\mathbb C$ in the lepton algebra.
By doing so, we lose $9$ (or perhaps even $10$) of the $18$ degrees of freedom, $8$ of which seem to re-appear in the odd part of the Dirac algebra.
It is very likely, therefore, that these $18$ dimensions contain
everything that we need for the rest of the bosonic part of the standard model. But the interpretation of these $18$
dimensions is likely to be different in some details from the standard model.

First let us define the elements $x,y,z$ of the group $Sym(4)$  to be the images of $i,j,k$, so that they satisfy the relations
$x^2=y^2=z^2=1$ and $xyz=1$.
These elements play the role of the three colours in the standard model, 
such that colour confinement is enforced by the condition $xyz=1$.
The standard model divides up the representation $3^++3^-$ on the basis of the eigenvalues of $\gamma_5$, that is $id$. On $3^+$ the $1$-eigenspace
has dimension $1$, and the $-1$-eigenspace has dimension $2$, while in $3^-$ it is the other way round. Therefore the standard model
breaks the symmetry $3+3$ down to $(1+2)+(2+1)$, which implies breaking the symmetry of $Sym(4)$ down to $Sym(2)\times Sym(2)$. 

The discrete model suggests that the symmetry might be restored by taking the up, down and strange quarks as a basis for $3^-$,
and the charm, bottom/beauty and top/truth quarks as a basis for $3^+$. In that way the representation $3^-$, which is the representation
of $Sym(4)$ as the rotation symmetry group of the cube, can perhaps be used to describe the symmetries of the baryon octet,
represented as vertices of the cube. The standard $2$-dimensional picture of the baryon octet \cite{GellMann}
based on the Weyl group
of $SU(3)$ is then simply the projection of the cube along one of the body diagonals. Of course, this copy of $SU(3)$ is not the same
as the gauge group $SU(3)$ that appears in the standard model, acting on the colours of the quarks.

The quarks themselves really only form a representation, not an algebra, 
and appear to use only $12$ of the $18$ dimensions, because colour confinement implies that there are
only two linearly independent colours. This reflects the fact that in the finite group, colour changes are three symmetries of order $2$.
Again this phenomenon arises from the use of the Lie group $SU(3)$ in place of the finite group $Sym(4)$, and the resulting
enforced symmetry-breaking down to $Sym(3)$. This is precisely analogous to the breaking of $Sym(3)$ down to $Sym(2)$
in the case of the weak interaction. Indeed, the same phenomenon can be seen in the necessary
breaking of $Sym(2)$ down to $Sym(1)$ for electromagnetism, which ensures that the finite group commutes with $U(1)$,
when $Sym(2)$ does not.

Therefore we have six spare dimensions, which are available in the finite case for elementary particles. 
I have already suggested that this part of the algebra might be used for the electromagnetic field, with a `Lorentz' gauge group
$SO(3,\mathbb C)$ in the standard model, or $SO(3)\times SO(3)$ in the group algebra model. In fact, the finite group contains
elements of determinant $-1$ in $3^+$, but not in $3^-$, so that we really need $SO(3)\times O(3)$, corresponding to
adjoining complex conjugation to $SO(3,\mathbb C)$. In other words, we have a charge-conjugation operator suitable for
building anti-particles.

But there is a difference between this and the Dirac algebra, that there is no Minkowski-space representation of $SO(3,1)$ here.
Therefore there is no mass here, and mass must be implemented somewhere else.
However, we can implement the momentum of the photons in the continuous action of $SO(3)$, 
as well as implementing polarisation in the finite group. 

\subsection{Mesons} 
If one (or both) of the representations $3^\pm$ is used for up, down and strange quarks, then the $Sym(4)$ symmetry attaches
signs to each of the quarks. These signs would appear to represent a discrete analogue of the standard model colours. They might
therefore be expected to play an important role in the modelling of composite particles such as mesons and baryons.

The meson octet of pions, kaons and eta meson is represented in the standard model via the adjoint representation of $SU(3)$,
made from the tensor product of the natural representation with its dual. We might therefore expect to find a corresponding representation
of the finite group by tensoring two $3$-dimensional representations. The available options are
\begin{eqnarray}
3^+\otimes 3^+ \cong 3^-\otimes 3^- &\cong& 3^- + (1^++2^0+3^+),\cr
3^+\otimes 3^-\cong 3^-\otimes 3^+ &\cong& 1^-+2^0+3^++3^-,
\end{eqnarray}
where the parentheses delimit the symmetric square. In all cases, therefore, the $8$-dimensional representation splits as $2^0+3^++3^-$.

Gell-Mann's original analysis of this case split the representation by restricting from $SU(3)$ 
to $SU(2)$, to give a structure $1+2+2+3$, allocated to the
eta meson, two pairs of kaons and three pions, in that order. 
The finite group representation splits differently, and would seem to require five kaons rather than four. An extra dimension for the kaons here
might indeed obviate the need for quantum superposition to explain the experimental observation of three distinct neutral kaon states.
However, this would require re-allocating the $\eta$ and $\eta'$ mesons, which may cause other problems. 
So this proposal can only be regarded as
very tentative at this stage.

\subsection{Baryons}
When it comes to baryons, the standard model uses 
the tensor cube,
which splits generically as $1+8+8+10$. 
Here $1$ denotes the anti-symmetric cube, $10$ denotes the symmetric cube and $8$ denotes what is sometimes called the `middle' cube.
Applying the same constructions to the finite group, we have
\begin{eqnarray}
\Lambda^3(3^+) &\cong& 1^-\cr
M^3(3^+) &\cong& 2^0+3^++3^-\cr
S^3(3^+)&\cong& 1^+ + 3^+ +3^++3^-,
\end{eqnarray}
and the same with all signs changed for $3^-$. 

In particular, the representation for the baryon octet is equivalent to that for the meson octet,
just as it is in the eightfold way. In the latter context, the octet divides into $2+3+3$ as follows:
\begin{eqnarray}
2 &\leftrightarrow& \Lambda, \Sigma^0\cr
3 &\leftrightarrow& p, \Sigma^-, \Xi^0\cr
3 &\leftrightarrow& n, \Sigma^+, \Xi^-
\end{eqnarray}
in such a way that the two triplets have the same total mass, as noted originally by Coleman and Glashow \cite{ColemanGlashow}.

Such a symmetry-breaking arises naturally in the finite group context. The triplet symmetries arise from cycling the quarks,
\begin{eqnarray}
u\rightarrow d\rightarrow s\rightarrow u. 
\end{eqnarray}
The doublet symmetry of $\Lambda$ and $\Sigma^0$ then appears to arise from swapping any two of the three quarks.
Note, incidentally, that the representation is not quite the same as the permutation representation on the vertices of the cube,
that I suggested earlier. The permutation representation has $1^++1^-$ in place of $2^0$.
The more subtle version using $2^0$ seems to be required for the physics here.

A similar picture emerges for the baryon decuplet, with
\begin{eqnarray}
1&\leftrightarrow& \Sigma^{*0}\cr
3&\leftrightarrow& \Delta^+,\Sigma^{*-},\Xi^{*0}\cr
3&\leftrightarrow& \Delta^0,\Sigma^{*+},\Xi^{*-}\cr
3&\leftrightarrow& \Delta^{++}, \Delta^-,\Omega^-
\end{eqnarray}
Again the two triplets with charges $0,+,-$ have the same total mass, to within experimental uncertainty.
The other triplet has a slightly smaller total mass, by about $(0.4\pm 0.1)\%$.
The baryon octet and decuplet between them use $18$ of the $27$ dimensions of the tensor cube, which begs the question whether the other 
$1+8=9$ dimensions can be interpreted in a similar way. 

\section{Towards a full model}
\label{fullmodel}
\subsection{The spinor algebra} 
\label{spinoralgebra}
There is no copy of the Lorentz group in the bosonic part of the group algebra that I have examined in detail above.
Therefore the theory so far is strictly non-relativistic. But there are several copies of the Lorentz group in the fermionic part of the algebra
\begin{eqnarray}
M_2(\mathbb C) + M_4(\mathbb R).
\end{eqnarray}
In particular, $M_2(\mathbb C)$ contains the Lorentz group in the form used by Dirac, that is $SL(2,\mathbb C)$, acting on the
Dirac spinor. There is another copy of $SL(2,\mathbb C)$ inside $M_4(\mathbb R)$, acting on the Dirac spinors
for the other two generations. And there is also a copy of $SO(3,1)$ inside $M_4(\mathbb R)$, 
which we can perhaps interpret as acting on the \emph{differences}
between the Dirac spinors for the three generations.

There are now three completely different groups that might in principle be interpreted as `the' Lorentz group, as well as
various ways in which they might be mixed together. It is not at all obvious which one to choose, and it is not even obvious whether there
is a viable choice for a single group that covers all of the functions of the Lorentz group in standard physics. 
It may rather be the case that different groups
can be taken as `the' Lorentz group in different contexts. 

One mathematical possibility, that may or may not have a
reasonable physical interpretation, is take a subgroup
\begin{eqnarray}
SL(2,\mathbb C) \times SO(3,1) & \subset & M_2(\mathbb C) + M_4(\mathbb R)
\end{eqnarray}
and then choose a homomorphism
\begin{eqnarray}
\varphi&:& SL(2,\mathbb C) \rightarrow SO(3,1)
\end{eqnarray}
so that there is a copy of $SL(2,\mathbb C)$ consisting of all elements of the form 
\begin{eqnarray}
(g,\varphi(g))&\in& SL(2,\mathbb C) \times SO(3,1).
\end{eqnarray}
A construction of this type would combine the quantum version $SL(2,\mathbb C)$ with the classical version $SO(3,1)$
in a way that would permit different observers 
to choose different
homomorphisms $\varphi$ appropriate to their particular circumstances. Or it may be that the homomorphism $\varphi$
can be fixed, and the different observers simply pick different copies of $SO(3,1)$ inside $M_4(\mathbb R)$.
The latter option could, at least in principle, provide a mechanism for incorporating general relativity into the theory, and allowing every
observer to choose their own definition of inertial frame, in their own gravitational context.

We can perhaps interpret $M_4(\mathbb R)$, in particle physics terms, as acting on the \emph{differences}
between the Dirac spinors for the three generations.
Since the differences between the three generations are mass differences, the latter interpretation implies that
$M_4(\mathbb R)$ acts on macroscopic mass, and therefore on momentum and energy. It follows that the representation
$4^0$ plays the role of quantising $4$-momentum and/or spacetime. One can therefore choose to gauge these parameters
with $SO(3,1)$, as is normally done in relativity, or with a larger group. The larger group $SL(4,\mathbb R)$ does not,
of course, preserve the (rest) mass as defined in special relativity, so that this model requires the rest masses to vary
with acceleration or the gravitational field.
But it must be stressed that $SL(4,\mathbb R)$ is a \emph{gauge group}, so that it doesn't change the theory, it only
changes the coordinates that are used to calculate in the theory. 

Certainly the spinor algebra bears some mathematical resemblance to the Dirac algebra, which is a copy of $M_4(\mathbb C)$.
The even part of the Dirac algebra consists of two copies of $M_2(\mathbb C)$, obtained by projections with $(1\pm \gamma_5)/2$.
The group algebra similarly contains two copies of $M_2(\mathbb C)$, one of them inside $M_4(\mathbb R)$. These can also be
obtained by projections, but of a rather different kind.

First of all, we project the group algebra with $(1-(-1))/2$ to obtain the fermionic part. Then to obtain $M_4(\mathbb R)$ we project with
\begin{eqnarray}\label{massidem}
&&(4 + w(1+i+j+k) + w^2(1-i-j-k))/6
\end{eqnarray}
and to obtain $M_2(\mathbb C)$ we project with
\begin{eqnarray}
&&(4 - 2w(1+i+j+k) - 2w^2(1-i-j-k))/6. 
\end{eqnarray}
The two possible choices of complex structure in the latter case are given by
\begin{eqnarray}
\mathrm i &=& \pm ( (j-k) + (i-j)w +(k-i)w^2) d/6\sqrt2.
\end{eqnarray}

Another significant difference between this spinor algebra and the Dirac algebra is that the odd part of the Dirac algebra
appears here only in a real form, not in a complex form. It is clear, therefore, that a direct translation between the two models
is not going to be possible. 
It may be, therefore, that the group algebra model can be ruled out on these grounds. But it may be that these differences can be
reconciled in some way, so it is worth pursuing the model a little further, rather than giving up at this stage.

\subsection{The Dirac equation}
Since the Dirac equation relates $4$-momentum to the Dirac spinors, it involves at least one representation other than
the $2^\pm$ that are used for the spinor. The equation must be invariant under the action of the finite group, and is therefore
a tensor product equation in the representation theory. Since both sides of the Dirac equation are spinors, we only need to
consider the tensors with the bosonic representations, that is $1^\pm$, $2^0$ and $3^\pm$: 

\begin{eqnarray}
2^+\otimes 1^+ \cong 2^-\otimes 1^-&\cong&2^+,\cr
2^+\otimes 1^- \cong 2^-\otimes 1^+&\cong&2^-,\cr
2^+\otimes 2^0 \cong 2^-\otimes 2^0 &\cong& 4^0,\cr
2^+\otimes 3^+ \cong2^-\otimes 3^- &\cong& 2^-+4^0,\cr
2^+\otimes 3^- \cong 2^-\otimes 3^+ &\cong& 2^++4^0.
\end{eqnarray}
Since $i,j,k$ act trivially in $1^\pm$ and $2^0$, the Dirac equation must involve at least one of the representations $3^\pm$.
One possible equation is then
\begin{eqnarray}
2^+\otimes 3^+ &=& 2^+\otimes(1^-+2^0),
\end{eqnarray}
but there are a number of other possibilities. Since the standard model unifies $1^\pm$ into a complex $1$-space,
and $3^\pm$ into a complex $3$-space, however, all these possibilities give essentially the same Dirac equation in the end.
In order to obtain the usual Dirac equation, the momentum is represented in $3^+$, and the energy in $1^-$, 
leaving $2^0$ as a `mass plane', for implementing various mass triplets that sum to a well-defined constant.
The two sides of the equation are not, however, Dirac spinors in the usual sense. To enforce this condition,
it is necessary to restrict from $G$ to $H$, so that $4^0$ can be identified
with the Dirac spinor $2b+2c$ for this subgroup. Thus we obtain three different Dirac equations
for the three different generations.

\subsection{Tensor products}
\label{tensors}
Other tensor products that might give rise to useful equations include
\begin{eqnarray}
\Lambda^2(2^0)&\cong&1^-\cr
S^2(2^0)&\cong& 1^++2^0\cr
2^0\otimes 3^+\cong 2^0\otimes 3^- &\cong & 3^++3^-\cr
4^0\otimes 2^0 &\cong& 2^++2^-+4^0\cr
4^0\otimes 3^+\cong 4^0\otimes 3^- &\cong & 2^++2^-+4^0+4^0\cr
2^+\otimes 4^0\cong 2^-\otimes 4^0 &\cong& 2^0+3^++3^-\cr
\Lambda^2(4^0) &\cong& 1^-+2^0+3^-\cr
S^2(4^0) &\cong & 1^++3^++3^++3^-
\end{eqnarray}
In particular, there is significant interplay between $2^0$ and $1^++1^-$ that seems worthy of further investigation, in addition to the
interplay between $3^+$ and $1^-+2^0$ already noted. The last two equations
suggest a possible alternative way of implementing a version of the Dirac algebra to exploit both of these. 

In fact, it is possible to write the entire group algebra as a tensor product, which may help to relate it to the tensor product structure
of the standard model. Looking first at the bosonic part, that is the group algebra of $Sym(4)$, we see a decomposition of the
representation as
\begin{eqnarray}
&2^0\otimes (1^\pm+2^0)\otimes (1^\pm+3^\pm)
\end{eqnarray}
in which the signs can be chosen arbitrarily. Combining the first and last terms gives a decomposition as
\begin{eqnarray}
&(1^\pm+2^0) \otimes (2^0+3^++3^-)
\end{eqnarray}
in which the two factors have already been identified as the finite equivalents of the adjoint representations of the gauge groups
$SU(2)$ and $SU(3)$ of the weak and strong interactions.

Turning now to the fermionic part of the algebra, there are many different decompositions, ranging from 
\begin{eqnarray}
&2^0\otimes 3^\pm \otimes 4^0,
\end{eqnarray}
with unbroken symmetry,
to 
\begin{eqnarray}
&(1^\pm+1^\pm)\otimes (1^\pm+2^0)\otimes (2^++2^-)
\end{eqnarray}
with completely broken symmetry. The symmetry of the individual factors can be broken independently, giving rise to
a total of eight different tensor decompositions. If we break the symmetry only of $3^\pm$, then we can combine this with
the bosonic part of the algebra into
\begin{eqnarray}
&2^0\otimes(1^\pm+2^0)\otimes (1^\pm+3^\pm+4^0).
\end{eqnarray}

This looks rather like three copies of a real version of the Dirac algebra, and there are various ways we might consider
reducing from the $6$ dimensions of 
\begin{eqnarray}
2^0\otimes(1^\pm+2^0) &=& 1^++1^-+2^0+2^0
\end{eqnarray}
down to two real or two complex dimensions. For example, we might consider $1^++1^-+2^0$ as a suitable way to
implement the action of the weak $SU(2)$ on left-handed but not right-handed spinors. The group algebra suggests, however,
that it is better to keep two copies of the left-handed spinors, in order to avoid the necessity for further symmetry-breaking.

Whatever choice we make here, the factor $1^\pm+3^\pm+4^0$ clearly divides the Dirac algebra into the odd part,
represented by $4^0$, and the even part, represented by $1+3$, in which the scalar and pseudoscalar ($1$ and $\gamma_5$)
are separated from the adjoint $SL(2,\mathbb C)$. In particular, there is an action of the Lorentz group $SO(3,1)$ on
$4^0$ which extends to an action of $GL(4,\mathbb R)$. Therefore an implementation of the Dirac algebra in this form,
if it can be carried out successfully, will be
not just Lorentz-covariant, but generally covariant.
In other words, this finite model of quantum mechanics may be consistent with general relativity, in a way that the continuous
standard model is not.

\subsection{Towards a mass gauge}
If so, then the matrix algebra $M_4(\mathbb R)$ describes all the mass/momentum/energy properties of elementary
particles, and the group $GL(4,\mathbb R)$ acts on it by left-multiplication to determine the
appropriate coordinates in any particular gravitational environment. The finite group $G$ acts on it by right-multiplication,
to permute the different types of elementary particles. Therefore the model permits four independent fundamental masses
(three fundamental mass ratios), from which all others can in principle be determined.

Such determination, however, requires us to identify the mass gauge group $GL(4,\mathbb R)$ with the group of
general covariance, without which there are $15$ or $16$ independent masses, as in the standard model.
With this identification, however, we should expect there to be four masses that, in the absence of a plausible quantum
theory of gravity, can only be determined by experiment, in such a way that all other masses can be obtained from
these four by calculations in the group algebra.
There is, of course, a wide choice of which four masses we consider to be fundamental. We seem to need a generation triplet
and a generation singlet, so my preferred choice is the three generations of electron, plus the proton. 
In addition, there is one very special mass,
corresponding to the identity matrix in $SL(4,\mathbb R)$.

This mass is therefore independent of the mass gauge, and hence can be
defined to be the same for all observers. The identity element is defined by the idempotent (\ref{massidem}). 
It is not at all clear
at this stage how to relate the terms in this formula to the particles, but there is certainly some generation symmetry, defined by
$1,w,w^2$, and some colour symmetry, defined by $i,j,k$. 
These symmetries are mixed together in a strange way, 
but there are six linearly independent colour terms, which seems to imply three protons, and three linearly
independent generation terms without colour, which seems to imply one electron from each generation.

Hence the conjecture must be that the sum of the masses of the three generations of electron, plus three protons, is the same
for all observers. Now it was observed in \cite{perspective} that this total mass is, experimentally, equal to the mass of five neutrons.
While this does not prove that this mass is invariant, it does add some credibility to the claim. Further development of the 
model
\cite{gl23model} provides additional justification.

More controversial is the claim that the individual masses of the electron and proton, relative to the neutron, are not invariant.
It must again be stressed that gauging these masses differently does not change the theory in any way, but merely changes
the coordinates. It is, of course, \emph{convenient} to use the mass coordinates that we measure in experiment, just as it is convenient to use
our conventional measurements of space and time in relativity, but it is not strictly \emph{necessary}. 

In order to analyse why and how we measure the particular mass values that we do, it is necessary in principle to analyse the whole history of physics,
and the decisions that were made at various points as to which measurements depend on which other measurements. These decisions are
regularly revised, and it is impossible in practice to disentangle the electron/proton/neutron mass ratios from a whole host of
other factors all over the theory of physics. Hence the best we can hope to do here is point out approximate
coincidences between these mass ratios and properties of the gravitational environment. Moreover, the accuracy we can expect cannot be
greater than that obtained in early direct experiments, rather than in modern indirect experiments that use a great deal more of the theory.
Such coincidences were pointed out in \cite{perspective}, to a relative accuracy of $10^{-5}$ for the proton/neutron ratio, and
$2\times 10^{-4}$ for the electron/neutron ratio.

\subsection{Chirality and helicity}
The terms helicity and chirality are generally used in physics to mean a handedness in $3$-dimensional and $4$-dimensional spacetimes respectively.
Mathematically, they are properties of orthogonal groups, usually $SO(2,1)$ and $SO(3,1)$ respectively. A finite group representation is always compact,
so we must work with $SO(3)$ and $SO(4)$ instead. Helicity can then be defined by the difference between the representations $3^+$ and $3^-$, since
in the latter case $G$ lies in $SO(3)$ and in the former case it does not.

Chirality then must be defined by the representation $4^0$, and the distinction between the two different $SO(3)$ quotients of $SO(4)$.
These quotients are best distinguished in the anti-symmetric square representation, which splits up as the sum of one $3$-dimensonal
representation of each quotient. In our case, we have $1^-+2^0+3^-$, so that one chiral piece is an irreducible representation $3^-$ of $Sym(4)$,
while the other is a reducible representation $1^-+2^0$ of $Sym(3)$. In other words, we have a `left-handed' chiral piece with broken
symmetry to represent the weak force, and a `right-handed' chiral piece with unbroken symmetry to represent part of electrodynamics
and/or the strong force. The latter is then combined with $3^+$ so that both helicities are included.

It is also possible to use the Dirac spinor in the form $2^++2^-$ for the same purpose. Here we find the left-handed part as $1^++1^-+1^-$ instead
of $1^-+2^0$. These two versions of chirality are the same on restriction to the subgroup $H$, which I have suggested as an analogue of the
generationless standard model. However, 
this implementation of chirality is completely different from the standard model
implementation, which is 
a distinction between $2b$ and $2c$, or between $2^+$ and $2^-$. 

Hence a direct translation from the
group algebra model to the standard model is impossible here, which may of course invalidate this proposed model.
On the other hand, chirality in the standard model is something that has to be added in by hand, using projections with $(1\pm \gamma_5)/2$,
whereas in the group algebra model it emerges naturally, complete with symmetry-breaking, so we should at least consider the possibility
that the new model is telling us something important at this point. 

The essential difference between the two models is in the treatment of spacetime.
In the standard model, spacetime is bosonic, whereas in the new model it is fermionic.
Now it is possible to regard bosonic spacetime as an obstacle to unification, in the following sense. With the standard 
implementation and interpretation of $SL(2,\mathbb C)$ as a double cover of the restricted Lorentz group $SO(3,1)^\circ$, an extension to
general relativity appears to require a double cover of $SL(4,\mathbb R)$ to act on the Dirac spinor. But there is no double cover
of $SL(4,\mathbb R)$ that has finite-dimensional representations, so a unification of this kind is impossible. 

A fermionic spacetime has no such problems, since in this case
both $SL(2,\mathbb C)$ and $SO(3,1)$ embed in $SL(4,\mathbb R)$ as subgroups. Of course, this breaks the traditional link between
$SL(2,\mathbb C)$ and $SO(3,1)$, so this may be a deal-breaker. But breaking this link is a necessary condition for an algebraic
unification of particle physics with general relativity, so maybe it is a price worth paying.
I have already discussed one potential resolution of this issue in Section~\ref{spinoralgebra}. 

\subsection{Quantization of spacetime}
Indeed, it could be said that the standard way of writing spacetime vectors $(x,y,z,t)$ as complex Hermitian matrices
\begin{eqnarray}
(x,y,z,t)&\mapsto&
\begin{pmatrix}
t+z & x+\mathrm iy\cr x-\mathrm iy & t-z
\end{pmatrix}
\end{eqnarray}
already does express spacetime fermionically. In this case the (bosonic) action of $SO(3,1)$ is expressed in matrix terms as
\begin{eqnarray}
X&\mapsto& MXM^\dagger
\end{eqnarray}
where $M^\dagger$ is the complex conjugate transpose of the matrix $M$ in $SL(2,\mathbb C)$.
A fermionic action is an action by multiplication, but multiplication by elements of $SL(2,\mathbb C)$ does not preserve the set of
Hermitian matrices, so has no meaning in the standard formalism.

But Hermitian matrices are in practice used in three different ways in physics. Given an Hermitian matrix $H$, one defines an
anti-Hermitian matrix
\begin{eqnarray}
A&:=&\mathrm iH
\end{eqnarray}
and a unitary matrix
\begin{eqnarray}
U&:=&\exp(A).
\end{eqnarray}
It is normal simply to write down $H$, rather than the implied matrices $A$ or $U$, and therefore not to distinguish too carefully between the
three. 

However, they belong to quite different mathematical structures: Hermitian matrices generate Jordan algebras, anti-Hermitian
matrices generate Lie algebras, and unitary matrices generate Lie groups. In the case of these spacetime matrices, 
the corresponding Lie group is $U(2)$, consisting of scalars $U(1)$ together with $SU(2)$. 
The Hermitian and anti-Hermitian matrices have no natural scale, so they do not quantise spacetime.
But the unitary matrices do have a natural scale, so in principle they could be used to quantise spacetime, via
\begin{eqnarray}
(x,y,z,t) &\mapsto& \begin{pmatrix} t+\mathrm iz & x+\mathrm iy\cr -x+\mathrm iy& t-\mathrm iz\end{pmatrix},
\end{eqnarray}
in such a way that the determinant is a natural Euclidean metric.

Now the action of $SL(2,\mathbb C)$
transfers 
from Hermitian matrices to anti-Hermitian matrices, but does not transfer to unitary matrices. For the action of
a  matrix $M$ to commute with exponentiation, it is necessary and sufficient that $M^\dagger M=1$, that is to say, $M$ is unitary.
In other words, we obtain an action of $U(2)$ by conjugation on itself, in place of the action of $SL(2,\mathbb C)$. These two actions
agree insofar as they both contain the same action of $SO(3)$ as an action of $SU(2)$ by conjugation. That is, they agree
as far as non-relativistic quantum mechanics is concerned.
 
Of course, the absence of an action of $SL(2,\mathbb C)$ on quantised spacetime may be a deal-breaker. But what we gain is
a multiplicative action of $U(2)$ on itself, such that $SU(2)$ can act by conjugation as the non-relativistic spin group, and by
multiplication as the weak gauge group. This overlap between the spin group and the gauge group may be what gives the weak force
its chirality.
 
 \section{Relations to other models}
\subsection{Clifford algebras and octonions}
The proposals made in Section~\ref{tensors}
to write the group algebra as a tensor product are rather different from the usual ideas involving Clifford algebras.
A Clifford algebra is essentially a tensor product of a (spin) representation with itself. It is clearly not sufficient in the group algebra
to take only a spin representation here, and we must have at least the sum of a bosonic and a fermionic representation.
Now it turns out that
\begin{eqnarray}
\Lambda^2(3^\pm+4^0) &=& 3^- + (1^-+2^0+3^-) + (2^++2^-+4^0+4^0)\cr
S^2(3^\pm+4^0) &=& (1^++2^0+3^+)+ (1^++3^++3^++3^-) \cr && + (2^++2^-+4^0+4^0)
\end{eqnarray}
so that the tensor square of $3^\pm+4^0$ consists of a copy of the group algebra plus the trivial representation $1^+$.

Of course, this only gives the structure as a representation of $G$, not the structure as an algebra, so that the latter needs to be
put back in somehow. Nevertheless, this structure contains a $3\times 3$ component which looks suitable for colours and/or generations,
two $3\times 4$ components that look suitable for (a) the $12$ fermion labels and (b) spinors, and a $4\times 4$ component suitable for
the Dirac algebra.

Without the guidance provided by the finite group $G$, however, it is hard to know what structure to put on this space. It is very tempting, for example,
to implement the representation $3^-+4^0$ as imaginary octonions, either complex, or split or compact real, as 
Dixon \cite{Dixon}, Furey \cite{Furey1} and others have done.
These and other approaches then extend from $7$ to $8$ dimensions in order to translate into a Clifford algebra \cite{Furey2,Furey3,Clifford}, 
usually some real or complex
form of $Cl(6)$, but sometimes a larger Clifford algebra, or a tensor product of Clifford algebras  \cite{Todorov}. 

Both $3^-$ and $4^0$ are real orthogonal representations, in which all elements of $G$ have determinant $1$.
In the absence of $G$, therefore, it is reasonable to give them symmetry groups $SO(3)$ and $SO(4)$, and thereby
embed $G$ in $SO(7)$.
We can then restrict to $SO(6)$, and re-introduce spinors using $Spin(6)\cong SU(4)$.
This is in effect what the Pati--Salam model \cite{PatiSalam} 
does, also adding in $SU(2)\times SU(2)$ to make up the total adjoint dimension of $21$. Thus
the whole group is $Spin(4)\times Spin(6)$ rather than $Spin(7)$.

Alternatively, we can embed $G$ in $SO(3,4)$ in order to incorporate some Lorentz transformations. 
At one time (unpublished lecture at EPFL, Lausanne, 7th June 2017) I advocated using this
group to attempt a unification of particle physics with relativity. This group can also be implemented in terms of the split real octonions if desired.
Most approaches would then lift to the spin group $Spin(3,4)$ in order to provide Dirac spinors for the model, but in my approach
the representation  $4^0$ already contains spinors, so it is not necessary to provide any more.

\subsection{General relativity}

The principal reason for making this suggestion was that on restriction from $SO(3,4)$ to $SO(3,3)$ the $7$-dimensional
representation breaks up as a scalar (mass) plus the anti-symmetric square of Minkowski space, and the symmetric square of the latter
consists of a scalar plus the Riemann Curvature Tensor \cite{thooft}. Moreover, to achieve this reduction to 
general relativity \cite{GR1,GR2}, it is necessary to
throw away the spinor representation $4^0$. Hence any unification along these lines must put the spinor back, in the form of
a vector ($7$-dimensional) representation of $SO(3,4)$, and extend the $20$-dimensional Riemann tensor to a $27$-dimensional
representation. This contrasts with general relativity itself, and the standard approaches \cite{GL4R1,GL4R2} to quantisation thereof,
which use $Spin(3,3)\cong SL(4,\mathbb R)$ instead of $SO(3,4)$.

If it is indeed possible to identify the Riemann curvature tensor as a representation of $G$ in this way, then it must be obtained by taking out
$1^++3^-+4^0$ from the symmetric square, so it appears as
\begin{eqnarray}
1^++2^0+3^++3^++3^++2^++2^-+4^0.
\end{eqnarray}
This can be arranged into a semblance of the traditional division into the Ricci scalar, plus Einstein tensor ($9$-dimensional real)
plus Weyl tensor ($5$-dimensional complex) as
\begin{eqnarray}
1^++(3^++3^++3^+)+(2^0+2^++2^-+4^0),
\end{eqnarray}
for example, although there are other possibilities, such as
\begin{eqnarray}
1^++(2^0+3^++4^0) + (2^++3^+)+(2^-+3^+).
\end{eqnarray}
Of course, we then have to combine this with $3^-+4^0$, and even then 
there is no guarantee that this leads to a viable quantisation of general relativity.

It may be significant that the first of these proposals for a finite version of the Einstein tensor, that is $3^++3^++3^+$, uses the representation
$3^+$, which lies inside $O(3)$ but \emph{not} inside $SO(3)$. This is the closest one can get to implementing the Lorentz group
$SO(3,1)$ using a finite group. Thus the Einstein tensor corresponds to one of the components $M_3(\mathbb R)$ of the group algebra.
In \cite{finite} I suggested using the corresponding part of the group algebra of the binary tetrahedral group for a quantisation of gravity,
using three generations of neutrinos as mediators.  Much the same can be done here, with the added advantage that there is a
second copy of $M_3(\mathbb R)$ available for the strong force.

The controversial use of fermions as mediators here appears to be forced on us by the fact that
the gauge group $GL(4,\mathbb R)$ of spacetime is fermionic in this model. 
It is also important to distinguish between $SO(3,1)$, extending to $SL(4,\mathbb R)$ or $GL(4,\mathbb R)$, which acts
on spacetime, and $GL(3,\mathbb R)$, which acts on the gravitational field. General relativity uses the same group for both purposes.
By separating them, we are able to distinguish between momentum and energy of massive fermions, described by $SO(3,1)$, extending to
$SL(4,\mathbb R)$, and
momentum of three generations of massless neutrinos, described by $GL(3,\mathbb R)$.

It may also be relevant that $SL(4,\mathbb R)$ factorises as a product of disjoint subgroups $SL(2,\mathbb C)$ and $GL(3,\mathbb R)$.
Transferring this result to the bosons gives a factorisation of $SO(3,3)$ as a product of $SO(3,\mathbb C)\cong
SO(3,1)$ and $GL(3,\mathbb R)$. This is not a
direct product, so that there is a lot of `mixing' between aspects of reality that are described by the group $SO(3,\mathbb C)$ and those that are
described by $GL(3,\mathbb R)$. Since these groups are not quantised by the finite group, they have relevance to macroscopic physics
rather than quantum physics. The group $SO(3,\mathbb C)$ is used in classical electrodynamics as the symmetry group of the
electromagnetic field. Hence $GL(3,\mathbb R)$ must surely act as the symmetry group of the gravitational field.

If so, then the mixing of these two groups implies that the electromagnetic field, quantised by photons, is affected by gravity, so that
light rays are bent by gravitational fields. Conversely, it implies that the gravitational field is affected by electromagnetism. The former effect
has been confirmed by direct observation, while the latter has not. Mathematically, indeed, one can choose a fixed copy of $GL(3,\mathbb R)$
for gravitation, so that only $SO(3,\mathbb C)$ varies according to the gravitational context. Or one can fix $SO(3,\mathbb C)$ and allow
$GL(3,\mathbb R)$ to vary. What one cannot do, and what physicists apparently want to do, is to fix both at the same time.

Indeed, I have already suggested that the electromagnetic field may be affected by gravity in more subtle ways than the bending of light,
in that the apparent anomaly in the muon magnetic moment might be an unrecognised gravitational effect. Conversely, is it possible that
observed gravitational anomalies, such as the flyby anomaly \cite{flyby} or the curious behaviour of galaxies \cite{Milgrom1}, 
might be related to 
an electromagnetic influence on the gravitational field?

\subsection{Quantum fields}
In order to 
address this question of the `mixing' of gravity with electromagnetism, we need to return to the group algebra,
and quantise the theory with the finite group. 
The basic operation that converts a particle into a field (or wave-function) is the tensor product with spacetime, represented by $4^0$.
Since this representation is fermionic, the tensor product operation maps fermionic particles to bosonic fields and vice versa.
Hence the particles and fields may not necessarily appear in quite the part of the algebra that one would expect.
We have already seen this effect in the case of the elementary fermions, in which the particles appear in the bosonic part of the
algebra, while the wave-functions appear in the fermionic part.

If we apply the same idea to photons, then they should appear as particles in the fermionic part of the algebra, and as waves
in the bosonic part. This puts them into $2^\pm$, with the sign possibly available for polarisation (as a discrete particle property).
Tensoring with $4^0$ gives $2^0+3^++3^-$, as an $8$-dimensional complex representation, or $16$-dimensional real. In other words,
there are two copies of $2^0+3^++3^-$ here, enough for a $6$-dimensional electromagnetic field, plus an $8$-dimensional gluon field,
and a $2$-dimensional mass-charge field, for example. 

Other interpretations are certainly possible, but this seems like a reasonable first
approximation. It at least allows all of the standard model bosonic gauge groups to appear and mix together
in a single representation. But it leaves no room for a bosonic theory of quantum gravity. Therefore this model only
permits a fermionic theory of quantum gravity. That is, the mediators are fermions, therefore presumably neutrinos,
and must live as particles in the bosonic part of the algebra, and as waves in the fermionic part. I have already
tentatively identified the neutrinos as living (as particles) somewhere in $1^++1^-+2^0+2^0$. 

Now tensoring $4^0$ with $1^+$ or $1^-$ just gives
back another copy of $4^0$, which seems unlikely to give us any physical process by which gravity can act.
Moreover, this part of the lepton algebra does not model the three generations, so is unlikely to model
mass, and therefore gravity, correctly.
So the most reasonable hypothesis is to put the neutrinos, as particles, into $2^0$, so that they appear as wave-functions
in two copies of $2^++2^-$. But $2^\pm$ are complex representations, and our tensor product is real, so that we do not
have two copies of the Dirac spinor that we might expect, but only two copies of a Majorana spinor.

On one hand, this is good, since it models neutrinos with only half as many degrees of freedom as the charged leptons.
On the other hand, it is bad, because it does not give us room for three linearly independent generations. But recall that $2^0$
is the representation as the symmetry group of an equilateral triangle. Therefore there is a perfectly good discrete symmetry
between three generations, even though these generations are not linearly independent. 

Now consider what happens when a neutrino propagates as a wave, and is detected as a particle in an experiment. At that moment of
detection, the wave-function in $2^++2^-$ `collapses' back into a part of the tensor product $2^0\otimes 4^0$. Here $4^0$ represents the
local quantised spacetime, and the detection picks out a particular vector in $2^0$ as the generation of the neutrino. There are two very
important things to note here. First, that there is no well-defined mathematical operation of `un-tensoring' a tensor product. Taking the
tensor product necessarily loses scalar information which it is impossible to recover. Second, that the collapse of the wave-function
depends on the local definition of spacetime.

There are various ways to interpret this result, depending on whether one regards the particle or the wave as being primary.
In the standard wave-function model, the problem is to get the scalar mass back out of the tensor product, which does not
contain that information. In the particle model, the neutrino has an intrinsic generation label in $2^0$, which we cannot determine
directly, but which we can probe by tensoring with our ambient spacetime $4^0$, 
in such a way that
what we get out of the model is a set of probabilities for the various possible outcomes. 

Whichever way we look at it, it is hard to
avoid the conclusion that the measured generation of a neutrino depends on the environment. 
At first glance, this seems like an absurd conclusion. It certainly is not consistent with the original version of the standard model,
in which neutrinos with zero mass had an intrinsic and immutable generation. But this version of the standard model was long ago
proved to be incorrect, as experiment has demonstrated conclusively that neutrinos generated in the sun as electron neutrinos are
detected on Earth in equal numbers in the three generations. So one of the two assumptions is wrong: either the neutrinos have
non-zero mass, or the generation label is not intrinsic. The 
standard model adopts the former hypothesis,
while my model adopts the latter.

\subsection{Wave-functions}
This description of photons and neutrinos suggests a more general treatment of wave-particle duality in terms of taking tensor
products with quantised spacetime in $4^0$. That is, from any particle lying in a representation $R$, we should
expect the corresponding wave-function
to lie in the representation $R\otimes 4^0$. 
In other words, $R$ contains the intrinsic properties of a particle, and $R\otimes 4^0$ contains the properties of its interactions
with its environment.

As a general principle, all this means is that the discrete properties of a particle `in isolation',
described by the representation $R$, become continuous properties of a particle in a given (macroscopic) spacetime environment,
described by the representation $R\otimes 4^0$. Admittedly, this is a very different interpretation from the usual textbook one, in which
the wave-function in $R\otimes 4^0$ is regarded as being intrinsic to the particle, and the `collapsed' wave-function in $R$
is regarded as a function of the experimental measurement. 
But it is not a new interpretation, 
and was suggested
by Rovelli \cite{RQM}, who calls it Relational Quantum Mechanics (RQM).

In particular, if the wave-function is a property of the embedding of the particle in its environment, rather than an intrinsic property
of the particle itself, then the measurement becomes a process of determining the intrinsic discrete value from the way the
particle interacts with its environment. The intrinsic values may be re-distributed amongst the elementary particles, but all are
discrete, and all are
conserved. The purpose of quantum mechanics, and particle physics more generally, is to try to predict which particles end up
with which quantum numbers.

Since we cannot isolate a single interaction from the rest of the universe, we cannot know what will happen in that interaction.
But by using the wave-function, we can work out what happens on average to a large number of elementary particles in
`similar' (defined by the experiment) circumstances. That is exactly what quantum mechanics does, and always has done,
very successfully. 

As many people have discovered over
nearly a century of grappling with the problem, it does not make sense to regard the wave-function as being intrinsic to the
elementary particle. Ultimately, in principle, the wave-function is a property of the entire universe.
(That is why the many-worlds interpretation of quantum mechanics is  logically unassailable, even if
 it is felt by many to be philosophically absurd.)
Nevertheless, the wave-function is almost universally regarded, in practice, as intrinsic to the elementary particle.
All we really have to do is abandon this patently absurd assumption, and then, as Rovelli points out,
to all intents and purposes, 
the measurement problem ceases to be a problem at all.

\section{Quantisation}
\subsection{Quantising the algebra}
The finite group $G$ quantises all the matrix algebra summands of the group algebra, by mapping to a finite set of matrices in
each algebra. To some extent this quantisation can be transferred from the algebra to the representations, but there is no
well-defined scale factor on any individual representation. This is essentially the same problem that occurs in trying to
`un-tensor' a tensor product, and is a familiar phenomenon in quantum mechanics. The most well-known example is
the Heisenberg Uncertainty Principle: neither position nor momentum can be individually quantised, but the duality between
them is quantised by Planck's constant.

Let us look first at the example $M_2(\mathbb R)$. The $6$ matrices can be taken as
\begin{eqnarray}
\begin{pmatrix}1&0\cr0&1\end{pmatrix}, &
\frac12\begin{pmatrix}-1&\sqrt 3\cr -\sqrt 3 &-1\end{pmatrix},\quad
\frac12\begin{pmatrix}-1&\sqrt 3\cr \sqrt 3 &1\end{pmatrix},\cr
\begin{pmatrix}1&0\cr0&-1\end{pmatrix},&
\frac12\begin{pmatrix}-1&-\sqrt 3\cr -\sqrt 3 &1\end{pmatrix},\quad
\frac12\begin{pmatrix}-1&-\sqrt 3\cr \sqrt 3 &-1\end{pmatrix}.
\end{eqnarray}
Now we can quantise the representation $2^0$ by the $1$-eigenvectors of the three reflections, say on the scale
\begin{eqnarray}
(2,0) & (-1 , \sqrt 3) & (-1, -\sqrt3)
\end{eqnarray}
or by the $-1$-eigenvectors, on a possibly different scale, 
\begin{eqnarray}
\lambda(0,2) & \lambda(\sqrt 3, -1) & \lambda(-\sqrt 3, -1). 
\end{eqnarray}
Therefore, we have two sets of three well-defined directions in the representation $2^0$, each of which can be interpreted in
a variety of different ways according to the context.

 It is also possible to express all of this formalism in terms of complex numbers.
If we define
\begin{eqnarray}
\mathrm w &:= &(-1+\mathrm i\sqrt 3)/2
\end{eqnarray}
then we can take the two sets of three directions as
\begin{eqnarray}
1,\mathrm w, \mathrm w^2,&&\mathrm i,\mathrm{iw},\mathrm{iw}^2,
\end{eqnarray}
and define the action of $w$ as multiplication by $\mathrm w$, and the action of $d$ as complex conjugation.
In this case it makes sense to scale one set of vectors by a factor of $\sqrt 3$, relative to the other set. For example,
if we take
\begin{eqnarray}
2/\sqrt 3,2\mathrm w/\sqrt3, 2\mathrm w^2/\sqrt3,&&2\mathrm i/3,2\mathrm{iw}/3,2\mathrm{iw}^2/3,
\end{eqnarray}
and interpret the imaginary part as the electric charge, then we have charges $0,1,-1$ in the first set, and $2/3,-1/3,-1/3$
in the second.

\subsection{Quantising spacetime}
Now let us look at $M_4(\mathbb R)$. Here we want to interpret two dual copies of $4^0$ as a Euclidean version of
spacetime, and its dual mass-momentum. This is not exactly the same as the standard Minkowski version, which cannot be
quantised by a finite group. 
The quaternionic description given earlier shows that we can take the generators $w$ and $jd$ as the matrices
\begin{eqnarray}
\begin{pmatrix}1&0&0&0\cr 0&0&1&0\cr 0&0&0&1\cr 0&1&0&0\end{pmatrix},
\begin{pmatrix}0&0&0&-1\cr 0&0&-1&0\cr -1&0&0&0\cr 0 & 1 & 0 & 0\end{pmatrix}. 
\end{eqnarray}
These matrices demonstrate that Euclidean spacetime can be effectively quantised with hypercubes, though the scale of the
hypercubes is not defined by the algebra. 

However, it is important to note that $G$ is not the full rotation group of the hypercube, which contains also the \emph{left}
multiplications by $i,j,k$, so has order $192$. The reflection group has order $384$, and is equal to the Weyl group of $SO(9)$.
The group $G$ is not a Weyl group of anything, although it is isomorphic to a double cover of the Weyl group of $SO(6)$.
It is absolutely crucial for the structure of quantum mechanics that we use exactly the group $G$, and nothing else. 
The extensions to a Weyl group of type $D_4$, $B_4$ or $F_4$, all of which are popular in Grand Unified Theories,
do not preserve the essential 
physical structures.

It is essential
for the chirality of the weak force, and for the spin $1/2$ properties of fermions, that the spins $i,j,k$ act on one side and not the other.
It is essential for the non-chirality of the strong force, and for the existence of three generations of fermions, 
that $w$ and $d$ act on both sides, by conjugation, and not by multiplication. Although it may be tempting to have both `left-handed'
and `right-handed' spins in the finite group, this is not consistent with experiment. Similarly, it may be tempting, but unwise, to separate $w$ and $d$ into
left and right multiplications, to get a group of order $2304$, which is an extension of the Weyl group of $F_4$ by an automorphism
of order $2$. Such `symmetries' do not preserve the hypercube, so do not preserve the quantum structure of spacetime.

As an aside, note that there is a representation of $G$ in which the action of $w$ is by right-multiplication instead of conjugation,
and all the other generators $i,j,d$ act as before. But this is the representation $2^\pm$, regarded as a real $4$-dimensional
representation, not $4^0$. As I have already remarked, the element $w$ (generation symmetry) 
is effectively not used in the standard model, so that these two representations can be regarded as equal, as long as we are 
not interested in mass. But if we are interested in generation symmetry and mass in general, then we must distinguish these
two representations carefully. In particular there is an invariant complex structure on $2^\pm$, unique up to complex conjugation, while there is
no invariant complex structure on $4^0$.  
With respect to the bases I have given, the complex structures are defined by 
taking one of the quaternions $\pm(j-k)/\sqrt 2$ as $\sqrt{-1}$. 

\subsection{Consequences}
Hence we see a natural chirality in this choice of complex structure. 
One might also interpret $w$ as being all that is left in the finite model of the weak
$SU(2)$ in the standard model, since $i,j,k$ are colours, and $d$ lies in $U(2)$ but not in $SU(2)$. If so, then we note that $w$ acts on
half of the representation $4^0$, fixing the other half, but it acts on all of the representation $2^\pm$. Hence the identification of
$4^0$ with $2^\pm$ creates an irreconcilable conflict with the underlying discrete structure of the three generations. 
The proposed finite model resolves this conflict by separating the $8$ real dimensions of the Dirac spinor into $2$ complex dimensions
and $4$ real dimensions, such that $w$ acts on $6$ of the $8$ real dimensions, rather than the $4$ in the standard model. It is this
extra triplet action on a $2$-space that enables the finite model to model the three generations.

Returning now to $4^0$, the four axes $\pm1$, $\pm i$, $\pm j$, $\pm k$ are determined by the group,
and because the scale is not defined, these axes can be identified in macroscopic Euclidean spacetime also, as defined by a 
particular observer or a particular experiment, as a suitable average over a sufficiently large number of elementary particles.
In order to interpret these directions, we need again to imagine `un-tensoring' the matrix algebra in the form $4^0\otimes 4^0$. In particular,
all absolute scale factors disappear, but \emph{relative} scale factors are still relevant. More importantly, the directions $i,j,k$ also become 
{only relative}, but relative to \emph{what} is not immediately clear.

For example, if we try to `untensor' $i,j,k$ as $1\otimes i$, $1\otimes j$, $1\otimes k$ then we might want to identify the right-hand factors as
momenta, in a suitable dimensionless form, and the left-hand factor therefore in the dual representation as a suitable dimensionless
version of time. The tensor product then represents a finite version of a differential equation, so that the individual terms are 
time derivatives of momenta, that is to say forces. Removing the mass scale, the forces become accelerations, which 
probably occur in two forms:
an acceleration due to gravity, and an acceleration due to rotation.

\subsection{Duality}
We have not yet removed the length scale, and in a certain sense, it cannot be removed. That is, there is a duality between the very small and
the very big, so that the smaller the quantum structures we want to investigate, the larger the macroscopic length scale must be.
Even for the study of the leptons, it would seem to be necessary to choose the directions $i,j,k$ to be the macroscopic directions
up/down, East/West and North/South, in some order and some orientation. Then it becomes impossible to study neutrino oscillations,
or the muon magnetic moment, or the decay of kaons into pions,
or anything else that involves particles outside the first generation,
without taking into account that these directions are not constant over the scale of the experiment.

For the study of quarks, or even for the study of the proton, it is surely necessary to go to a larger scale, probably the scale of the Solar System,
in order to see all the effects of accelerations due to the gravitational attractions of the Sun and the Moon, not to mention Jupiter
and the other planets. The weak force, for example, shows a symmetry-breaking from $3$ dimensions to $2+1$, which looks to be 
approximately dual to the
symmetry-breaking of space on the Solar System scale down to the $2$-dimensional ecliptic plane.
It then becomes impossible to study the electron/proton/neutron mass ratios, or the fine-structure constant, 
or the masses of the up and down quarks, or the masses of any other particles,
without taking into account the motion of the experiment relative to the Sun and the Moon.

All this discussion is aimed at trying to decide how to use the mass gauge algebra $M_4(\mathbb R)$ to gauge the masses of the
elementary particles. Since this gauge algebra is also the gauge for general relativity, it is impossible to separate the two aspects
of this algebra. Mathematically, there is only one possible conclusion: if there is a discrete model of quantum mechanics,
that is consistent with \emph{any} realistic theory of gravity that has an $SL(4,\mathbb R)$ gauge,
then
the masses of the elementary particles are determined by the large-scale structure of the gravitational field.

Of course, this is a highly controversial suggestion, and it is obvious that one cannot change any of the masses individually,
and one cannot change a collection of masses without changing many other parameters at the same time. But,
ultimately, there is only one equation that relates gravity to mass: 
therefore, if gravity is determined by mass, then mass is determined by gravity. These are simply two different ways
of giving a physical interpretation to the same mathematical equation.

\subsection{An alternative interpretation}
It is apparent that there are two different `levels' of Lie groups in the group algebra. At the top level are the 
unimodular subgroups (that is, with determinant of modulus $1$, or equal to $1$, if appropriate)
 of the full matrix algebras that
the group algebra decomposes into, and 
at  the lower level are the Lie groups that are actually generated by the finite group.
Let us tabulate them here for comparison:
\begin{eqnarray}
&\begin{array}{ccccc}
\mbox{Name} & \mbox{Algebra} & \mbox{Unimodular} & \mbox{Compact}&\mbox{Finite}\cr\hline
1^+&\mathbb R &SO(1) & SO(1) & Sym(1)\cr
1^-&\mathbb R &O(1) & O(1)& Sym(2)\cr
2^0&M_2(\mathbb R) & SL^\pm(2,\mathbb R) & O(2)&Sym(3)\cr
3^+&M_3(\mathbb R) &SL^\pm(3,\mathbb R) & O(3)&Sym(4)\cr
3^-&M_3(\mathbb R) &SL(3,\mathbb R) & SO(3)&Sym(4)\cr
2^\pm &  M_2(\mathbb C) &U(1)\circ SL(2,\mathbb C) & U(2)&G\cr
4^0 & M_4(\mathbb R) & SL(4,\mathbb R) & SO(4)&G\cr\hline
\end{array}
\end{eqnarray}

For quantum mechanics, we really only need the finite group, but the standard theory uses at least the compact group.
For purposes of quantum superpositions and measurements by macroscopic apparatus, we need to add together lots of
group elements, and therefore we need the full matrix algebra. The standard theory here removes the overall positive real
scale factor,
so that for this purpose it uses the unimodular group.  However, quantum superposition is always done using complex scalars
rather than real scalars, so that in some cases $SL(n,\mathbb R)$ turns into $SU(n)$.

In a similar way, the standard model uses the unimodular group $U(1)\circ SL(2,\mathbb C)$, acting on the Dirac spinor $2^++2^-$,
to model the fermionic properties of elementary particles that would appear from our mathematical analysis to 
require the compact group $U(2)\circ SO(4)$ acting on $2^\pm+4^0$. 
To achieve this feat, it first
converts to the bosonic form $SO(2)\times SO(3,1)$, and then 
identifies $SO(2)$ with $U(1)$ for complex numbers to 
change the signature of $SO(3,1)$ to $SO(4)$. 
But it does not extend $SO(2)$ to $SO(3)$, or $U(1)$ to $SU(2)$,
which is what would be necessary to model the three generations of fermions.

But the real problem here is that the model has converted the fermionic representation of $SO(4)$ to a bosonic
representation of $SO(3,1)$. This would seem to be a fatal error from which there is no possible recovery. This group $SO(4)$ is
required for quantisation of gravity, and is both compact and fermionic. 
To convert it into a non-compact group completely prevents quantisation of any kind.
And to convert it to a bosonic 
group converts the obvious force-carriers, which are neutrinos, into spin $2$ gravitons,
for which there is no experimental evidence whatsoever.
\section{Conclusion}
The group $G$ considered in this paper is the only finite group that lies in $U(2)$ but not in $SU(2)$, and also acts irreducibly in $SO(3)$.
Both properties seem to be physically important: the former is needed for electroweak unification, and the latter is required for
generation symmetry of quarks. Hence $G$ is the only plausible candidate for a finite group on which to base a discrete model of quantum mechanics.
Of course, there is no guarantee that such a model exists. 

In this paper, I have shown that most of the mathematical structures of the standard model of particle physics can be found
inside the real group algebra of $G$, but in many cases there are added twists. It can certainly be argued that these added twists
therefore invalidate the proposed model, but I have tried to argue that, on the contrary, these added twists are, at least in some
cases, improvements to the standard model. The most obvious of these is that the three generations are incorporated
naturally into the model, in such a way that the experimentally detected generation of a neutrino
can change without a corresponding change in mass. 
The chirality and symmetry-breaking of the weak force also are direct consequences of the
representation theory of the finite group.

\end{document}